\documentclass{article}
\usepackage[centertags]{amsmath}
\usepackage{amsthm,amsfonts}
\usepackage{amssymb}
\usepackage{latexsym}
\usepackage{color}
\usepackage{tabularx}
\usepackage{graphicx}
\usepackage{enumerate}
\usepackage{epstopdf}
\usepackage{tikz}
\usepackage{booktabs}
\usetikzlibrary{arrows,positioning}
\definecolor{boxcolor}{HTML}{B9DCFF}
\usepackage{booktabs}
\usepackage[flushleft]{threeparttable}
\usepackage{caption}
\usepackage{pgfplots}
\usepackage[colorlinks=true]{hyperref}
\usepackage{cleveref}
\usepackage[position=t,singlelinecheck=off]{subfig}
\usepackage{subfig}
\usepackage{array, makecell}
\usepackage{indentfirst}

\makeatletter
\newcommand*{\rom}[1]{\expandafter\@slowromancap\romannumeral #1@}
\makeatother
\theoremstyle{plain}

\numberwithin{equation}{section}

\title{Studying the mixed transmission in a community with age heterogeneity: COVID-19 as a case study}
\author{\small{Xiaoying Wang $^{1*}$, Qing Han$^2$ \& Jude Dzevela Kong$^2$}\\
\small{$^1$ Department of Mathematics, Trent University}\\
\small{Peterborough, ON K9L 0G2, Canada}\\
\small{$^2$ Africa-Canada Artificial Intelligence and Data Innovation Consortium (ACADIC)}\\
\small{Laboratory for Industrial and Applied Mathematics (LIAM)}\\
\small{Department of Mathematics and Statistics, York University}\\
\small{Toronto, ON M3J 1P3, Canada}\\
\small{$^*$ Corresponding author: xiaoyingwang@trentu.ca}
}
\date{}

\begin{document}

\maketitle

\begin{abstract}
COVID-19 has been prevalent worldwide for about 2 years now and has brought unprecedented challenges to our society. Before vaccines were available, the main disease intervention strategies were non-pharmaceutical. Starting December 2020, in Ontario, Canada, vaccines were approved for administering to vulnerable individuals and gradually expanded to all individuals above the age of 12. As the vaccine coverage reached a satisfactory level among the eligible population, normal social activities resumed and schools reopened starting September 2021. However, when schools reopen for in-person learning, children under the age of 12 are unvaccinated and are at higher risks of contracting the virus. We propose an age-stratified model based on the age and vaccine eligibility of the individuals. We fit our model to the data in Ontario, Canada and obtain a good fitting result. The results show that a relaxed between-group contact rate may trigger future epidemic waves more easily than an increased within-group contact rate. An increasing mixed contact rate of the older group quickly amplifies the daily incidence numbers for both groups whereas an increasing mixed contact rate of the younger group mainly leads to future waves in the younger group alone. The results indicate the importance of accelerating vaccine rollout for younger individuals in mitigating
disease spread.
\end{abstract}


\vskip0.2cm
\baselineskip=18pt

\section{Introduction}
COVID-19 has been prevalent worldwide for around 2 years now since the initial identification in Wuhan, China in December 2019. The global pandemic has caused more than 328,000,000 total infections and over 5,500,000 deaths worldwide \cite{WHO}. In Ontario, Canada alone, the total confirmed case number of COVID-19 surpassed 970,000 and more than 10,000 population deceased due to complications of COVID-19 infection \cite{JHU}.

The novel pneumonia disease can be transmitted via close contact between susceptible and infected populations, similar to other pneumonia diseases such as influenza. However, recent evidences show that the majority of COVID-19 transmission may be attributed to the aerosol droplets \cite{A_P, W_C}. The transmission via aerosol implies that the spread of the disease may occur in a long-range, which signifies the difficulties of mitigating the disease spread.

In Ontario, before vaccines were approved, the main disease mitigation strategies were non-pharmaceutical, such as closing non-essential businesses, practising social distancing, setting limitations of indoor gathering size, requiring mandatory face-covering etc. The non-pharmaceutical strategies are effective in mitigating the disease spread but disrupt normal social activities and cause huge economic loss. Therefore, aggressive strategies can only be imposed in a short period and are not practical in long term. Progressing into the second year of the pandemic, vaccines were approved by Health Canada in December 2020 and quickly rolled out across the province.

In the initial stage of the vaccine rollout program, vaccines were administered to the senior population or other vulnerable individuals with underlying health conditions. Starting in late May 2021, Health Canada expanded the criteria so that all individuals above the age of 12 were eligible for receiving vaccines. The vaccine for children between the age of 5 to 11 was approved recently in late November 2021 but the rollout takes time and only $3.2\%$ of the children in the age group are fully vaccinated up to now \cite{ONdata}.

Clinical trail evidences show that vaccines offer a high protection efficacy for vaccinated individuals that reduces the probability of infection to a large extent \cite{WHO}. As the vaccine coverage steadily increased and reached a satisfactory level among eligible individuals, the province lifted majority of the restrictive measures in July 2021. During this time, unvaccinated children were on the summer break and weren't on high risks of transmitting the disease because of limited contacts within the household.

Starting in September 2021, schools reopen and in-person learning resumes. In-person learning environment inevitably creates a larger social gathering size, larger contact numbers between the children, and a higher infection probability in the indoor setting. Evidences show that children are at a lower risk of developing serious symptoms after contracting the virus and therefore the within-group transmission may not raise serious concerns. However, children may transmit the disease to elder family members within the household after school. Infected elder individuals are more likely to develop serious symptoms and need medical attention, which causes stress on the provincial health system.

Mathematical models have been used extensively in studying COVID-19 since the initial emergence of the disease, see \cite{He_S, Iwata_K, Liu_Z, Liu_Z2, M_I, Ruan_SG, Tang_B, W_X, Zhao_H} for example. The aforementioned studies focus on modeling the initial wave of the pandemic and provide valuable insights in short-term predictions and disease intervention strategies. Age-stratified models were proposed to study the heterogeneity in social contact patterns among a susceptible population in early stage of the disease \cite{D_N, F_R, G_V, M_Z}. Statistical analyses were applied in studying various aspects of the disease spread, such as identifying the impact of undiagnosed cases in COVID-19 \cite{H_A}, the importance of household transmission \cite{P_L}, and suggestions for testing policy \cite{L_J}. As the disease evolves, particularly with the rollout of vaccines and resumption of normal social activities, models need to be revisited to gain better understanding of the disease spread.

In this paper, we propose an age-stratified model where individuals in the older group are eligible for receiving vaccines while the individuals in the younger group are not. We fit the model to the COVID-19 data in Ontario, Canada for illustration. In Ontario, the schools reopen starting in September 2021 and the vaccines for children haven't been approved by then. The individuals above the age of 12 are eligible for receiving vaccines but breakthrough infections may occur due to relaxed social distancing. The main objective of this paper is to identify the relative importance of within-group transmission and between-group transmission and shed light on future disease mitigation strategies.

\section{Methods}

\subsection{Data}
\begin{figure}[htbp]
	\begin{center}
		\includegraphics[height=8.5cm,width=14.5cm]{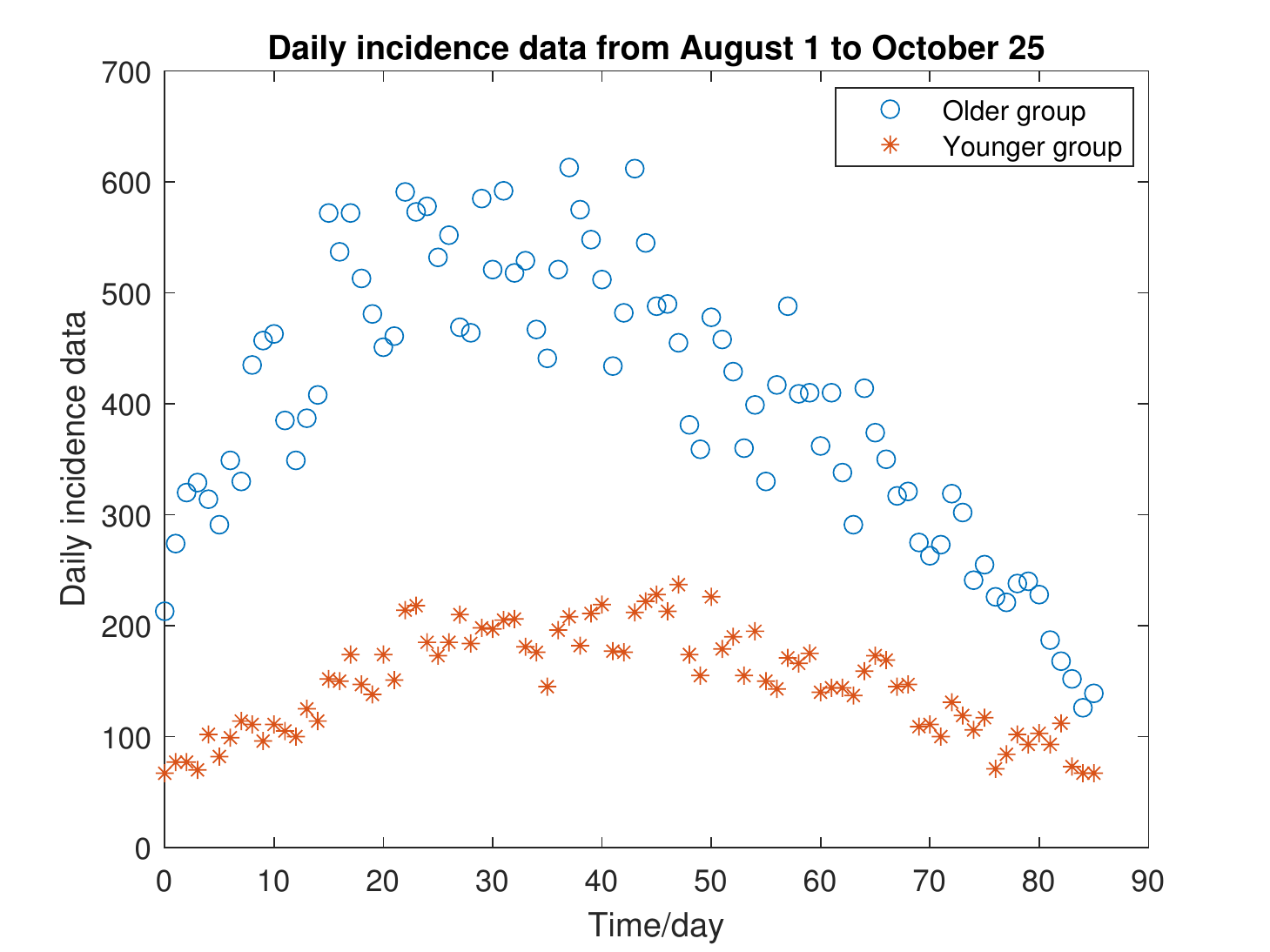}
		\caption{The daily incidence number for older/younger group in Ontario between August 1, 2021 and October 25, 2021.}
		\label{fig:p1}
	\end{center}
\end{figure}

We accessed the COVID-19 data in Ontario from the online data catalogue in Ontario government \cite{Data}. The data contains demographic information of individuals who are confirmed of positive COVID-19 infection among all 34 health units across Ontario, such as age, sex, and location of the reporting health unit. \autoref{fig:p1} shows the daily incidence number between Aug. 1, 2021 and Oct. 25, 2021, based on the estimated symptom onset date recorded in \cite{Data}. In the aforementioned time window, individuals who are above the age of 12 are eligible for receiving vaccines but individuals who are below the age of 12 haven't been approved for receiving vaccines yet. Therefore, in this study, we divide the entire population in Ontario into two groups: the older group and the younger group, where individuals are above/below the age of 12 respectively. \autoref{fig:p1} demonstrates the daily incidence case number for each group separately.

%P1

\subsection{The model}
We propose a compartmental model that stratifies the entire population into the following compartments: susceptible class ($S$), vaccinated class ($V$), exposed class ($E$), asymptomatic class ($A$), unreported symptomatic class ($U$), confirmed infected class ($I$), recovered class ($R$), and the deceased class ($D$). Moreover, within each compartment, we further stratify individuals into the older group and the younger group, as denoted by $S_i$ for $i=1,2$ for example. Figure \ref{fig:model} demonstrates the flowchart of the model. Note that within the time window between August and October, the individuals below the age of 12 are not eligible for receiving vaccines and therefore, the compartment $V$ comprises of vaccinated individuals of the older group only.

%P2
\tikzset{block/.style={draw, fill=boxcolor, minimum size=2em, text width =
1.5cm, align = center, minimum height = 1.5cm},
arrow/.style={->, black, text = black, line width=0.3mm}}
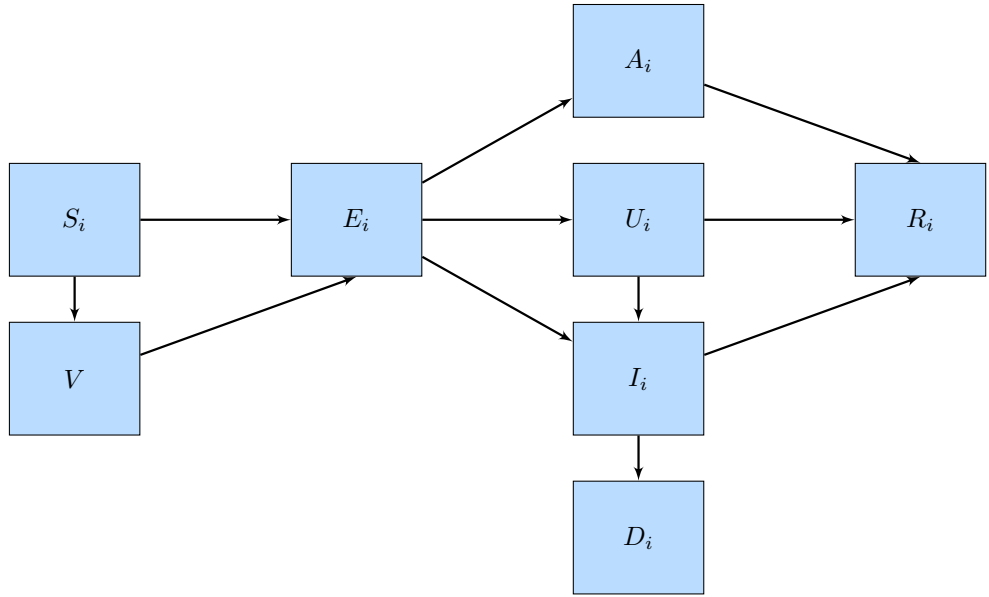
\begin{figure}[htbp]
\centering
\begin{tikzpicture}[auto, >=latex']
    \node [block] (s) {$S_i$};
    \node [block, right = 2.0cm of s] (e) {$E_i$};
    \node [block, above right = 0.6cm and 2.0cm of e] (a) {$A_i$};
    \node [block, below right = 0.6cm and 2.0cm of e] (i) {$I_i$};
    \node [block, below right = 0.6cm and 2.0cm of a] (r) {$R_i$};
    \node [block, right = 2.0cm of e] (u) {$U_i$};
    \node [block, below = 0.6cm of i] (d) {$D_i$};
    \node [block, below = 0.6cm of s] (v) {$V$};

    \draw [arrow] (s) -- node [above] {} (e);
    \draw [arrow] (e) -- node [above] {} (a);
    \draw [arrow] (e) -- node [below] {} (i);
    \draw [arrow] (e) -- node [above] {} (u);
    \draw [arrow] (a)-- node [above] {} (r.north);
    \draw [arrow] (u) -- node [above] {} (r);
    \draw [arrow] (i) -- node [below] {} (r.south);
    \draw [arrow] (u) -- node [below] {} (i);
    \draw [arrow] (i) -- node [below] {} (d);
    \draw [arrow] (s) --node  [below] {} (v);
    \draw [arrow] (v) --node  [below] {} (e.south);
\end{tikzpicture}
\caption{The flow chart for model.}
\label{fig:model}
\end{figure}

In the model, individuals in the exposed class ($E$) are individuals who have been infected with the COVID-19 virus but are still in the incubation period and are not contagious yet. Vaccinated individuals are at a lower risk of contracting the virus but may still become infected because the vaccine efficacy is not $100\%.$ Asymptomatic individuals in class ($A$) are those who are contagious but show no symptoms of infection. On the other hand, individuals in class ($U$) show mild symptoms and are also contagious but do not seek medical resources for diagnosis. Individuals in class ($I$) are those who are confirmed with the COVID-19 infection and are in quarantine at home if they experience mild to moderate symptoms or are in the hospital if they show severe symptoms. Because individuals in class ($I$) follow the public health protocols, we assume that they maintain a very low possibility of transmitting the disease. Individuals in all the asymptomatic classes, unreported symptomatic class, and reported infected class may recover and are no longer transmissible and therefore move to the recovered class ($R$).

Based on the aforementioned discussions, the model for the older group, where individuals are above the age of 12 and are eligible to be vaccinated is
\begingroup
\allowdisplaybreaks
\begin{align}\label{modeli}
&\frac{dS_1}{dt}=-\Lambda_{11}-\Lambda_{12}-\omega S_1, \nonumber\\
&\frac{dV}{dt}=\omega S_1-\Lambda_{11}^v-\Lambda_{12}^v, \nonumber \\
&\frac{dE_1}{dt}=\Lambda_{11}+\Lambda_{12}+\Lambda_{11}^v+\Lambda_{12}^v-\mu_1 E_1, \nonumber \\
&\frac{dA_1}{dt}=\mu_1\rho_1E_1-\eta_{11}A_1,\\
&\frac{dU_1}{dt}=\mu_1\rho_2E_1-\eta_{12}U_1-\tau_1U_1, \nonumber \\
&\frac{dI_1}{dt}=\mu_1\left(1-\rho_1-\rho_2\right)E_1+\tau_1 U_1-\eta_{13}I_1-\delta_1 I_1, \nonumber \\
&\frac{dR_1}{dt}=\eta_{11}A_1+\eta_{12}U_1+\eta_{13}I_1, \nonumber \\
&\frac{dD_1}{dt}=\delta_1 I_1. \nonumber
\end{align}
\endgroup
In \eqref{modeli}, $\Lambda_{11}=\beta_{11} A_1 S_1+\beta_{12} U_1S_1$ represents the within-group transmission of the susceptible population, where $\beta_{11}, \beta_{12}$ are the infectious contact rates between susceptible individuals and asymptomatic/unreported symptomatic individuals in the older group. Moreover, $\Lambda_{12}=\beta_{13}A_2S_1+\beta_{14}U_2S_1$ represents the between-group transmission, where $\beta_{13}, \beta_{14}$ are the contact rates between susceptible individuals in the older group and the asymptomatic/unreported infected individuals in the younger group.

Similarly, $\Lambda_{11}^v=\beta_{11}^v(1-\epsilon)A_1V+\beta_{12}^v(1-\epsilon)U_1V$ is the within-group transmission of the vaccinated individuals, where $\beta_{11}^v, \beta_{12}^v$ are the contact rates between vaccinated individuals and asymptomatic/unreported symptomatic individuals in the older group. Furthermore, $\Lambda_{12}^v=\beta_{13}^v(1-\epsilon)A_2V+\beta_{14}^v(1-\epsilon)U_2V$ is the between-group transmission of the vaccinated individuals, where $\beta_{13}^v,\beta_{14}^v$ are the contact rates between the vaccinated individuals and the asymptomatic/unreported symptomatic individuals in the younger group. Because vaccinated individuals are at lower risk of transmitting the disease, we introduce the factor $(1-\epsilon)$ in $\Lambda_{11}^v,\Lambda_{12}^v$ to indicate the reduced probability of infection where $\epsilon$ is the vaccine efficacy and $0<\epsilon<1.$

The model for the younger group is similar but without the vaccinated compartment because individuals below the age of 12 are not eligible for receiving vaccines. The model is
\begin{equation}\label{modelii}
\aligned
&\frac{dS_2}{dt}=-\Lambda_{21}-\Lambda_{22},\\
&\frac{dE_2}{dt}=\Lambda_{21}+\Lambda_{22}-\mu_2 E_2,\\
&\frac{dA_2}{dt}=\mu_2\rho_1E_2-\eta_{21}A_2,\\
&\frac{dU_2}{dt}=\mu_2\rho_2E_2-\eta_{22}U_2-\tau_2U_2,\\
&\frac{dI_2}{dt}=\mu_2\left(1-\rho_1-\rho_2\right)E_2+\tau_2 U_2-\eta_{23}I_2-\delta_2 I_2,\\
&\frac{dR_2}{dt}=\eta_{21}A_2+\eta_{22}U_2+\eta_{23}I_2,\\
&\frac{dD_2}{dt}=\delta_2 I_2,
\endaligned
\end{equation}
where $\Lambda_{21}=\beta_{21} A_2 S_2+\beta_{22} U_2S_2$ is the within-group transmission of the susceptible individuals, and
$\Lambda_{22}=\beta_{23} A_1S_2+\beta_{24} U_1S_2$ is the between-group transmission of the susceptible individuals.

In \eqref{modeli}-\eqref{modelii}, $\rho_1$ is the proportion of the exposed individuals who move to the asymptomatic class and $\rho_2$ is the proportion of the exposed individuals who move to the unreported symptomatic class. This leaves $1-\rho_1-\rho_2$ as the proportion of the exposed individuals who eventually move to the reported symptomatic class. We assume that $\rho_1, \rho_2$ are the same for both the older and the younger groups because there is no evidence so far that suggests the proportions of asymptomatic individuals and symptomatic individuals differ significantly among different age groups.

Susceptible individuals in the older group transfer to the vaccinated class at a rate of $\omega.$ The parameter $\omega$ will be estimated below from the data in \autoref{fig:p1}. Exposed individuals in the older and the younger group move to either the asymptomatic class/unreported symptomatic class/reported infected class at a rate of $\mu_i$ for $i=1,2$ respectively. Based on \cite{He_X, Li_Q}, we fix $\mu_i=1/3$ for $i=1,2$ because evidence shows that individuals who have contracted the virus generally have an incubation period of 2-3 days before they become contagious.

Asymptomatic individuals in the older group and the younger group move to the recovered class at a rate of $\eta_{11}$ and $\eta_{21}$ respectively. Similarly, the unreported infected individuals and the confirmed infected individuals in the older/younger group recover from the infection at a rate of $\eta_{i2}$ and $\eta_{i3}$ for $i=1,2$. We assume that $\eta_{i,j}=1/7$ where $i=1,2, j=1,2,3$ based on \cite{Liu_Z}, which indicates that in average, infected individuals may transmit the disease within a period of 7 days until they recover and are no longer contagious.

Unreported infected individuals may seek diagnosis from medical facilities if the symptoms persist or their heath deteriorates over time and therefore move to the reported infected class at a rate of $\tau.$ In \cite{Li_Q}, the authors analyzed the patients' data and recorded that infected individuals waited for an average of 4.6 days after observing symptoms before they sought diagnosis. Hence, in our study, we assume that $\tau_i=1/4.6$ for $i=1,2.$

\subsection{The reproduction number}

We first calculate the basic reproduction number $\mathcal{R}_0$ for the age-structured model \eqref{modeli}-\eqref{modelii}. Direct calculations show that
\begin{equation*}
\frac{d(S_1+V+E_1+A_1+U_1+I_1+R_1+D_1)}{dt}=0,
\end{equation*}
which leads to the total population of the older group as a constant $N_1.$ Similarly, the total population in the younger group also remains as a constant $N_2.$ Hence, the disease-free equilibrium of \eqref{modeli}-\eqref{modelii} is
\begin{equation*}
\aligned
& \left(S_1,V,E_1,A_1,U_1,I_1,R_1,D_1,S_2,E_2,A_2,U_2,I_2,R_2,D_2\right)\\
&=\left(0,N_1,0,0,0,0,0,0,N_2,0,0,0,0,0,0\right).
\endaligned
\end{equation*}
Following the next generation matrix method for compartmental models in \cite{VD_P}, we obtain
\begin{equation*}
\aligned
& F_1=\Lambda_{11}+\Lambda_{12}+\Lambda_{11}^v+\Lambda_{12}^v,\quad F_2=0, \quad F_3=0, \quad F_4=0,\quad F_5=\Lambda_{21}+\Lambda_{22},\\
& F_6=0,\quad F_7=0,\quad F_8=0,\\
&V_1=\mu_1 E_1,\quad V_2=\eta_{11}A_1-\mu_1\rho_1E_1,\quad V_3=\eta_{12}U_1+\tau_1U_1-\mu_1\rho_2E_1,\\
&V_4=\eta_{13}I_1+\delta_1I_1-\mu_1(1-\rho_1-\rho_2)E_1-\tau_1U_1,\quad V_5=\mu_2E_2,\quad V_6=\eta_{21}A_2-\mu_2\rho_1E_2,\\
&V_7=\eta_{22}U_2+\tau_2U_2-\mu_2\rho_2E_2,\quad V_8=\eta_{23}I_2+\delta_2I_2-\mu_2(1-\rho_1-\rho_2)E_2-\tau_2U_2.
\endaligned
\end{equation*}
It follows that
\begin{equation}\label{r0}
\mathcal{R}_0=\frac{J_{11}+J_{55}+\sqrt{(J_{11}+J_{55})^2-4(J_{11}J_{55}-J_{15}J_{51})}}{2},
\end{equation}
where
\begin{equation*}
\aligned
& J_{11}=\frac{\beta_{11}^v(1-\epsilon)\rho_1N_1}{\eta_{11}}+\frac{\beta_{12}^v(1-\epsilon)\rho_2N_1}{\eta_{12}+\tau_1},\quad
J_{15}=\frac{\beta_{13}^v(1-\epsilon)\rho_1N_1}{\eta_{21}}+\frac{\beta_{14}^v(1-\epsilon)\rho_2N_1}{\eta_{22}+\tau_2},\\
& J_{51}=\frac{\beta_{23}\rho_1N_2}{\eta_{11}}+\frac{\beta_{24}\rho_2N_2}{\eta_{12}+\tau_1},\quad
J_{55}=\frac{\beta_{21}\rho_1N_2}{\eta_{21}}+\frac{\beta_{22}\rho_2N_2}{\eta_{22}+\tau_2}.
\endaligned
\end{equation*}

To predict the severity of the disease spread at an arbitrary time of the epidemics, we calculate the effective reproduction number
$\mathcal{R}_t=(S_t/S_0)\mathcal{R}_0$, where $S_t$ is the total susceptible population at time t and $S_0$ is the total susceptible population at the initial time when the epidemic starts \cite{Diekmann, Nishiura_H1}. Biologically, the effective reproduction number
implies that the disease persists if $R_t>1$ and diminishes if $R_t<1.$

\subsection{Parameter estimation}

Based on the census data in Ontario, the total population of the older group at the beginning of the pandemic is 12,932,471 and the total population of the younger group is 1,801,543. The number of cumulative infections until July 31, 2021 for the older group is 555,583 whereas the total infection case for the younger group is 50,156 \cite{ONdata}. The vaccine coverage data shows that a total number of 10,582,731 individuals received at least one dose and a total of 9,147,534 individuals are fully vaccinated with two doses. We average the data to obtain the estimated initial population of the vaccinated compartment $V_0=9,865,132.$ It follows that the initial susceptible population for the older group is $S_1(0)=2,511,756.$ Similarly, the initial susceptible population for the younger group is $S_2(0)=1,751,387.$ On August 1, the reported infection number for the older group is 200 and the reported infection number for the younger group is 18, which leads to $I_1(0)=200$ and $I_2(0)=18$ respectively.

We fit the model \eqref{modeli}, \eqref{modelii} to the daily incidence data in \autoref{fig:p1} by the Markov Chain Monte Carlo (MCMC) method and adopt the adaptive Metropolis-Hastings algorithm to carry out this approach \cite{H_H}. We run the algorithm for 10,000 iterations with a burn-in of the first 7,000 iterations. Geweke convergence test is employed to diagnose the convergence of the Markov chains. The estimated parameters and the initial data are listed in \autoref{table:para}.

\begin{table}[ht!]
\caption{Parameter estimates for the COVID-19 epidemics in Ontario, Canada}
\centering\renewcommand\cellalign{c}
\scalebox{0.8}{
\begin{tabular}{c c c c c}
\hline\hline
Parameter &  Definition & \makecell{Estimated \\ Mean Value} & \makecell{Standard \\ Deviation} & \makecell{Data \\ Source} \\
\hline
$\beta_{11}$ & \makecell{Contact rate between\\ $S_1$ and $A_1$} & $2.1772 \times 10^{-8}$ & $1.708 \times 10^{-9}$ & Fitted \\
$\beta_{13}$ & \makecell{Contact rate between\\ $S_1$ and $A_2$} &$4.8079 \times 10^{-7}$ & $3.5895 \times 10^{-8}$ & Fitted \\
$\omega$ & \makecell{Transition rate from \\ $S_1$ to $V$} & $3.7286 \times 10^{-2}$ & $4.3604 \times 10^{-3}$ & Fitted\\
$\beta_{11}^v$ & \makecell{Contact rate between\\ $V$ and $A_1$} & $5.5607 \times 10^{-7}$ & $2.4585 \times 10^{-8}$ & Fitted \\
$\beta_{13}^v$ & \makecell{Contact rate between\\ $V$ and $A_2$} & $1.4986 \times 10^{-8}$ & $1.1856 \times 10^{-9}$ & Fitted  \\
$\epsilon$ & \makecell{Contact rate reduction\\ between $V$ and $A_i/U_i$} & $0.96783 $ & $4.423 \times 10^{-3} $ & Fitted  \\
$\mu_1$ & \makecell{Transition rate from \\ $E_1$ to $A_1/U_1/I_1$} & $1/3$ & $-$ & \cite{Li_Q}\\
$\mu_2$ & \makecell{Transition rate from \\ $E_2$ to $A_2/U_2/I_2$} & $1/3$ & $-$ & \cite{Li_Q}\\
$\rho_1$ & \makecell{Proportionality of \\ transferred $E_i$ to $A_i$}  & $0.3$ & $-$ & \cite{Nishiura_H}\\
$\rho_2$ & \makecell{Proportionality of \\ transferred $E_i$ to $U_i$}  & $0.3$ & $-$ & \cite{Nishiura_H} \\
$\eta_{i1}$ & Recovery rate of $A_i$ & $1/7$ & $-$ & \cite{Liu_Z}\\
$\eta_{i2}$ & Recovery rate of $U_i$ & $1/7$ & $-$ & \cite{Liu_Z}\\
$\eta_{i3}$ & Recovery rate of $I_i$ & $1/7$ & $-$ & \cite{Liu_Z}\\
$\tau_i$   & \makecell{Transition rate from \\ $U_i$ to $I_i$} & $0.21739$ &$-$ & \cite{Liu_Z} \\
$\delta_1$ & Disease death rate of $I_1$ & $ 1.0928 \times 10^{-5} $& $6.1912 \times 10^{-7}$  & Fitted \\
$\beta_{21}$ & \makecell{Contact rate between\\ $S_2$ and $A_2$} & $4.079 \times 10^{-9}$ & $3.765 \times 10^{-10}$ & Fitted \\
$\beta_{23}$ & \makecell{Contact rate between\\ $S_2$ and $A_1$} &$7.0164 \times 10^{-8}$ & $2.5621 \times 10^{-9}$ & Fitted \\
$\delta_2$ & Disease death rate of $I_2$ & $ 2.9753 \times 10^{-6} $& $4.1384 \times 10^{-7}$  & Fitted \\
\hline
Initial Value & Definition & Estimated Mean Value & Standard Deviation & Data Source\\
\hline
$S_1(0)$ & Initial susceptible population of group-1& $2.5117 \times 10^{6}$   & $-$ & \cite{ONdata}\\
$V(0)$ & Initial vaccinated population & $9.8651 \times 10^{6} $ & $-$ & \cite{ONdata}\\
$E_1(0)$ & Initial exposed population of group-1& 912.8571 & $-$ & Fitted\\
$A_1(0)$ & Initial asymptomatic population of group-1& $678.83$ & $63.497$ &Fitted\\
$U_1(0)$ & Initial unreported population of group-1& $42.791$& $4.3679$ & Fitted\\
$I_1(0)$ & Initial reported case number of group-1 & $200$ & $-$& \cite{ONdata}\\
$S_2(0)$ & Initial susceptible population of group-2 & $1.7513 \times 10^6$ & $-$ & \cite{ONdata}\\
$E_2(0)$ & Initial exposed population of group-2& $287.1429$ & $-$ & Fitted\\
$A_2(0)$ & Initial asymptomatic population of group-2& $165.88$ & $15.296$ &Fitted\\
$U_2(0)$ & Initial unreported population of group-2& $184.81$& $32.996$ & Fitted\\
$I_2(0)$ & Initial reported case number of group-2& $18$ & $-$& \cite{ONdata}\\
\hline
\end{tabular}
}
\label{table:para}
\end{table}

\section{Results}
\noindent
(\rom{1}) Fitting results

\begin{figure}[!htpb]
\captionsetup[subfigure]{labelformat=parens,labelfont=bf}
\centering
\subfloat[]{\includegraphics[width=0.5\textwidth]{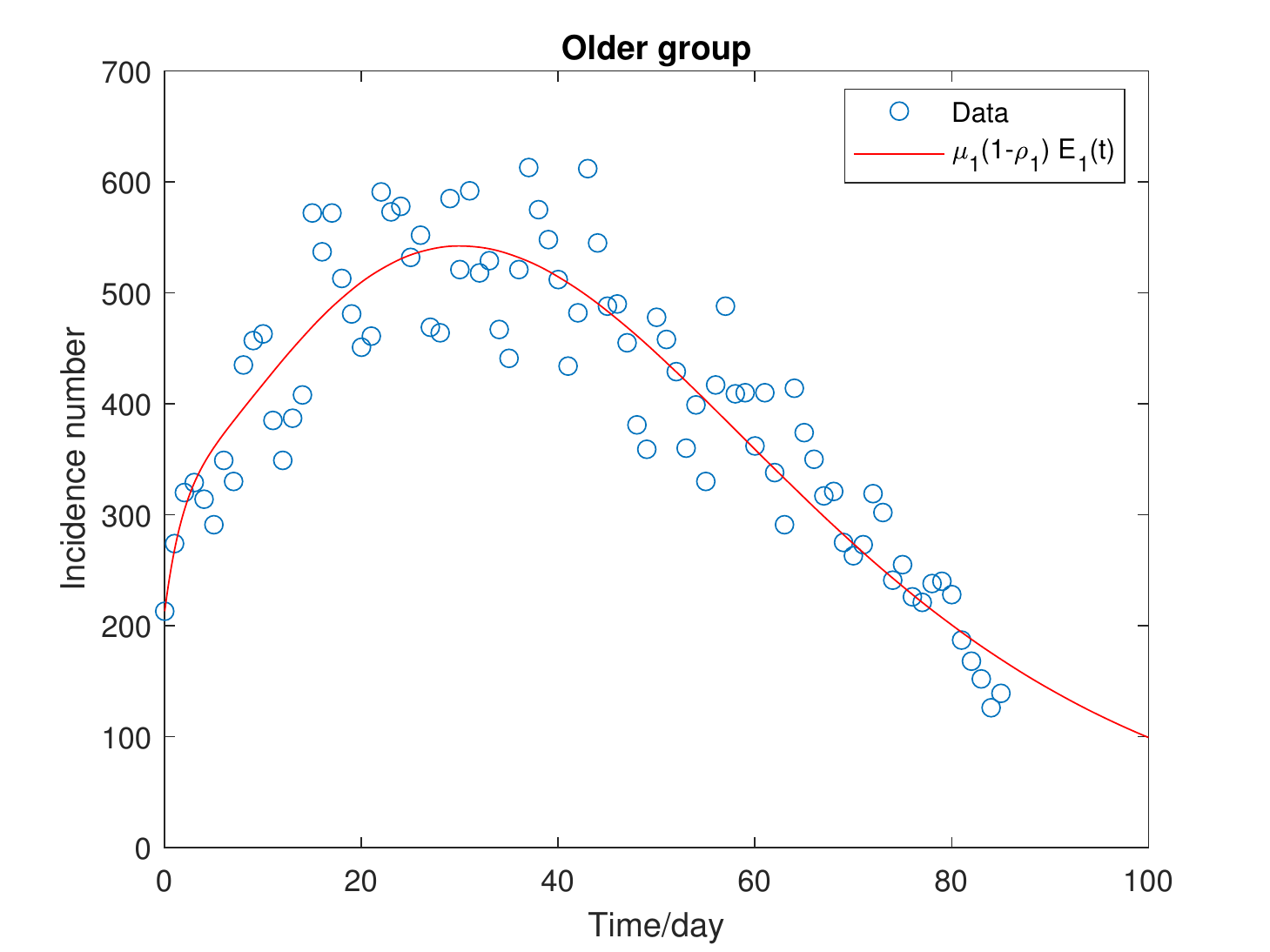}\label{fig:subfigure21}}\hfill
\subfloat[]{\includegraphics[width=0.5\textwidth]{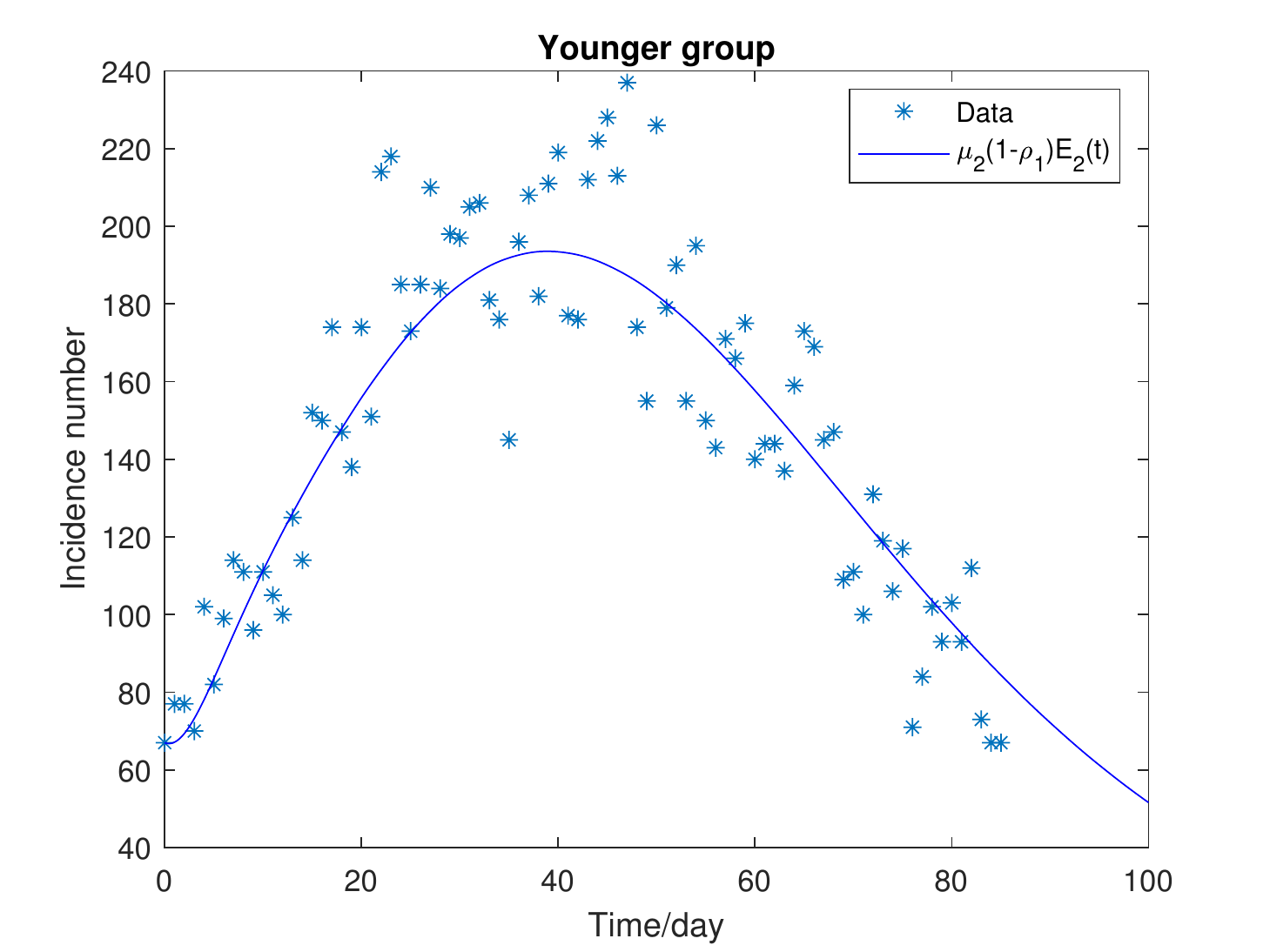}\label{fig:subfigure22}}\\
 \caption{The fitting result of model \eqref{modeli}-\eqref{modelii} to the daily incidence data between August 1 and October 25 in Ontario, Canada. \autoref{fig:subfigure21} shows the fitting of the daily incidence data to the older group. \autoref{fig:subfigure22} shows the fitting of the daily incidence data to the younger group.}\label{fig:p2}
\end{figure}

\autoref{fig:p2} shows the fitting result of the daily incidence data in \autoref{fig:p1} to the older group \eqref{modeli} and the younger group \eqref{modelii} respectively. As shown in \autoref{fig:p2}, the daily incidence data increases initially but then declines. The data shows the fourth wave of COVID-19 in Ontario, which is before the new variant emerges in Canada. The initial increase of the daily incidence case number may be attributed to the reopening in August and September, when the vaccine program steadily rolls out. It is then followed by a gradual decline, which is due to a combination of a few factors, such as a high vaccine coverage among eligible individuals, implementation of health protocols etc. \autoref{fig:p2} indicates that the daily incidence data changes in the same pattern for both the older and the younger groups but the infected case number peaks slightly later in the younger group than the older group.

Under the current set of parameters in \autoref{table:para}, the prediction of \eqref{modeli}, \eqref{modelii} indicates that the disease dies out eventually. The prediction is confirmed by the effective reproduction number $\mathcal{R}_t=0.31727<1$ by substituting the parameter values in \autoref{table:para}. However, starting in November and afterwards, vaccinated individuals may experience higher probabilities of infection because of the gradual decline of vaccine protection efficacy. Meanwhile, before December, individuals in the younger group are still not eligible for receiving vaccines but are at higher risks of infection. This is mainly because children stay indoors more often due to the weather condition and will contract the virus more easily if they follow the same health protocols as before. We will explore different scenarios below by increasing the within-group contact rate and the between-group contact rate.

%P3

\noindent
(\rom{2}) Comparing the relative importance of the within-group and between-group transmission

Next, we analyze the impact of $\beta_{11}$ and $\beta_{13}$ on future epidemic waves and compare the relative importance of the within-group contact rate and the between-group contact rate on the final epidemic size. By increasing the contact rate $\beta_{11}$ or $\beta_{13}$, we obtain a similar pattern of the future waves for both the older group and the younger group and hence show only the result of the older group.

\begin{figure}[!htpb]
\captionsetup[subfigure]{labelformat=parens,labelfont=bf}
\centering
\subfloat[]{\includegraphics[width=0.5\textwidth]{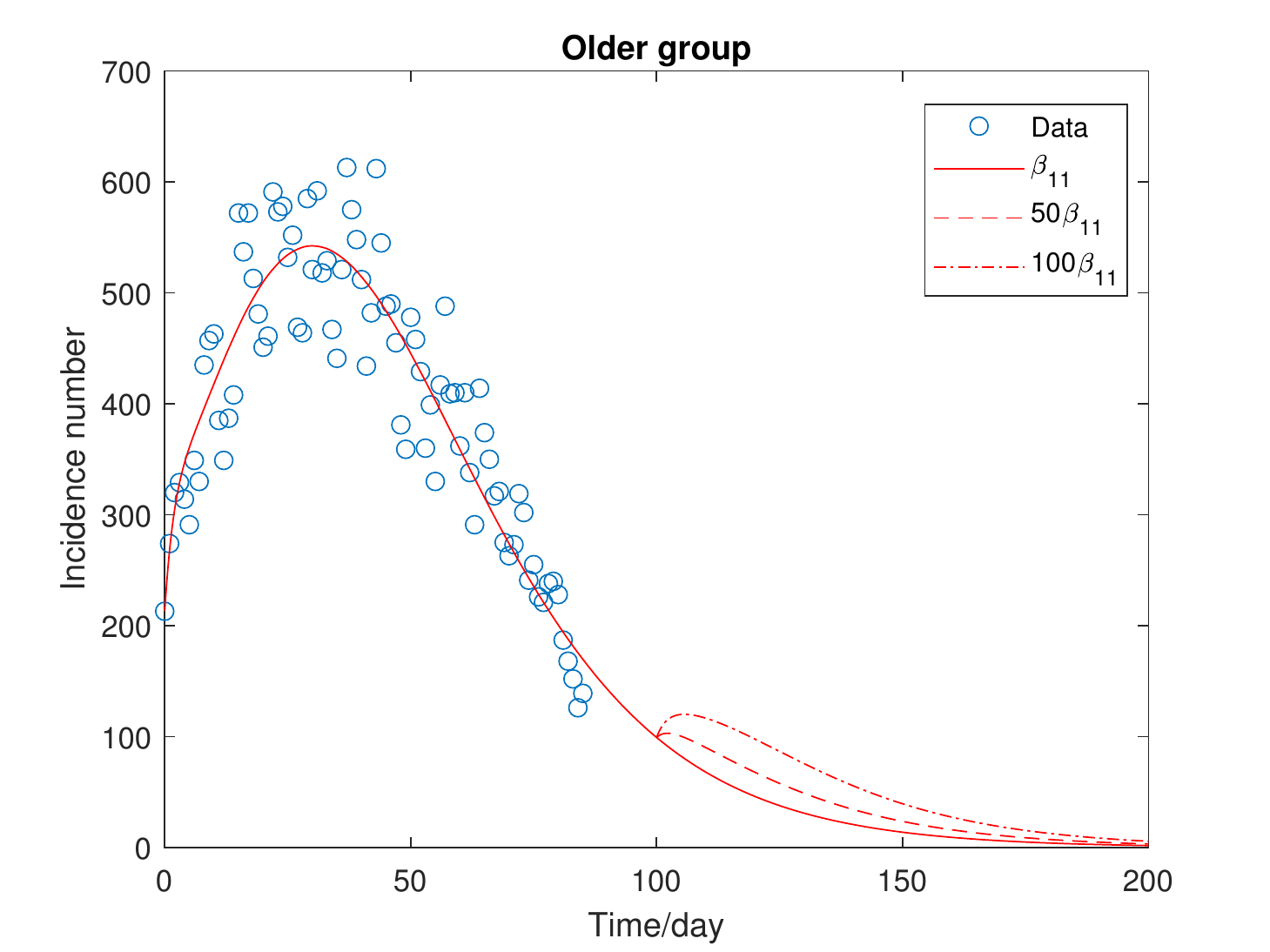}\label{fig:subfigure31}}\hfill
\subfloat[]{\includegraphics[width=0.5\textwidth]{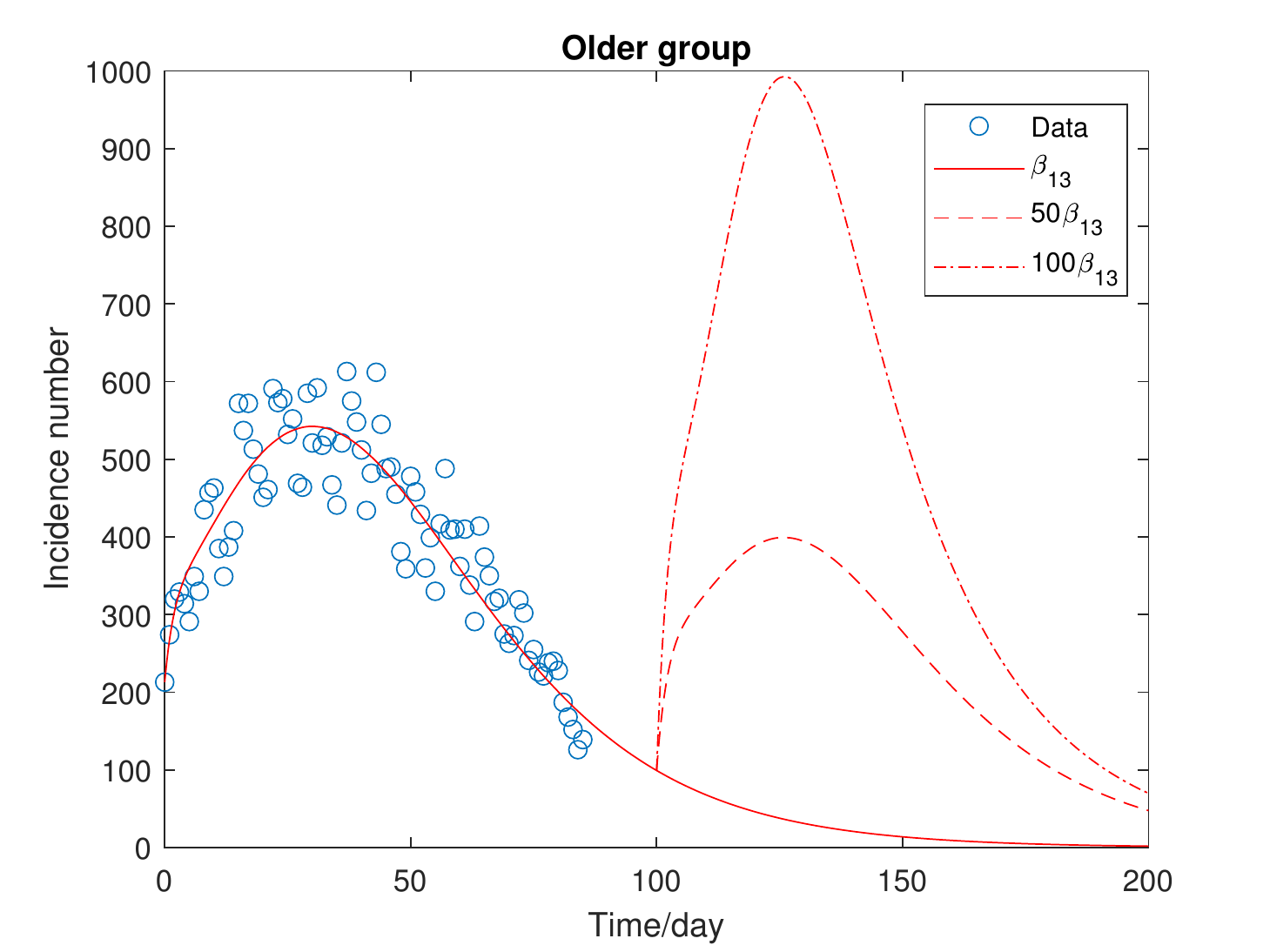}\label{fig:subfigure32}}\\
 \caption{Daily incidence number of the older group for different contact rate $\beta_{11}$ and $\beta_{13}$ respectively. All the parameters remain unchanged as shown in \autoref{table:para} except $\beta_{11}$ and $\beta_{13}$ respectively.}\label{fig:p3}
\end{figure}

\autoref{fig:subfigure31} shows that the daily incidence number of the older group increases and initiates a further wave if the within-group contact rate $\beta_{11}$ increases. The further epidemic wave achieves a higher peak if the contact rate $\beta_{11}$ is larger. However, the peak size of epidemics is relatively small compared to the peak size of the fourth wave even if the contact rate is larger and is on an order of different magnitude.

However, \autoref{fig:subfigure32} shows that an increasing between-group contact rate $\beta_{13}$ triggers a future epidemic wave of a large size. The peak size of future epidemic waves possibly exceeds the peak size of the fourth wave if the between-group contact rate of the older group is sufficiently large. \autoref{fig:p3} shows that the between-group contact rate of the older group imposes a much larger impact on future epidemic waves compared to the within-group contact rate.

\begin{figure}[!htpb]
\captionsetup[subfigure]{labelformat=parens,labelfont=bf}
\centering
\subfloat[]{\includegraphics[width=0.5\textwidth]{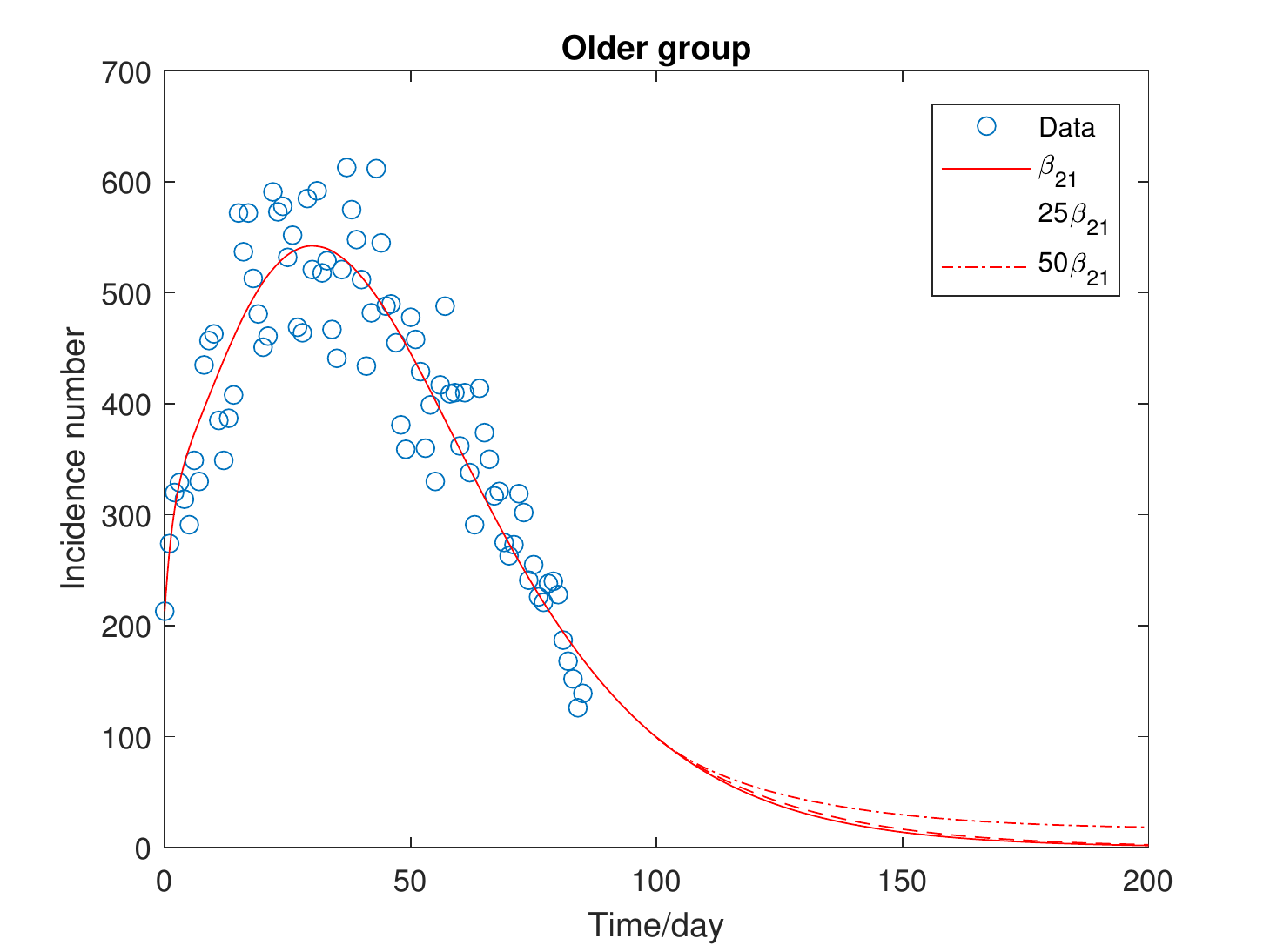}\label{fig:subfigure41}}\hfill
\subfloat[]{\includegraphics[width=0.5\textwidth]{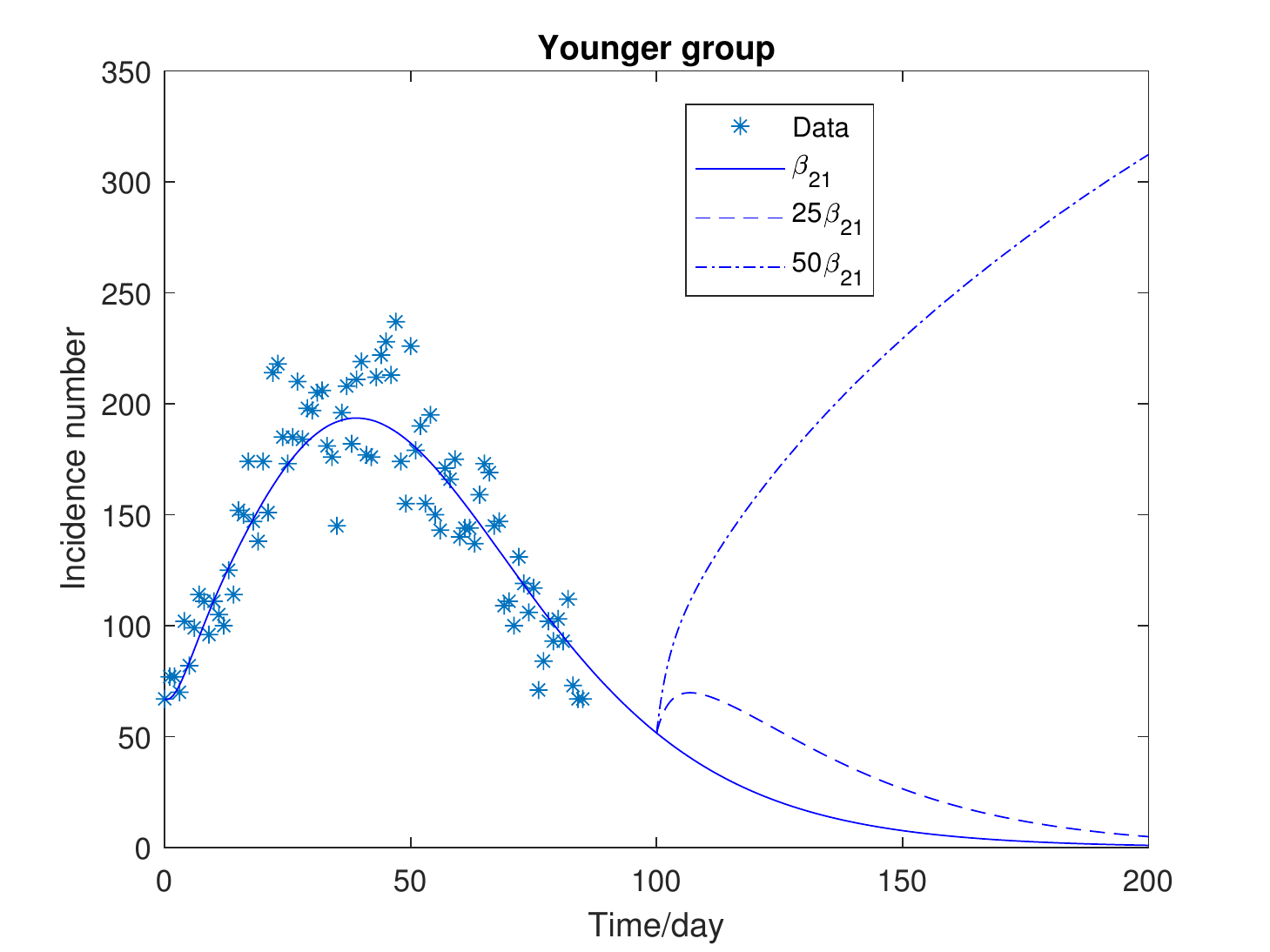}\label{fig:subfigure42}}\\
 \caption{Daily incidence number of the older group and the younger group respectively for different contact rate $\beta_{21}$. All the parameters remain unchanged as shown in \autoref{table:para} except $\beta_{21}.$}\label{fig:p4}
\end{figure}

Next, we increase the within-group contact rate of the younger group $\beta_{21}.$ \autoref{fig:subfigure41} demonstrates that the daily incidence number of the older group declines monotonically even if $\beta_{21}$ increases significantly. However, the daily incidence number of the younger group shifts to a gradual increase from the original decline if $\beta_{21}$ increases, as shown in Figure \autoref{fig:subfigure42}. More importantly, the daily incidence number of the younger group increases rapidly and forms a new epidemic wave of the size much larger than the fourth wave if $\beta_{21}$ is relatively large.

\begin{figure}[!htpb]
\captionsetup[subfigure]{labelformat=parens,labelfont=bf}
\centering
\subfloat[]{\includegraphics[width=0.5\textwidth]{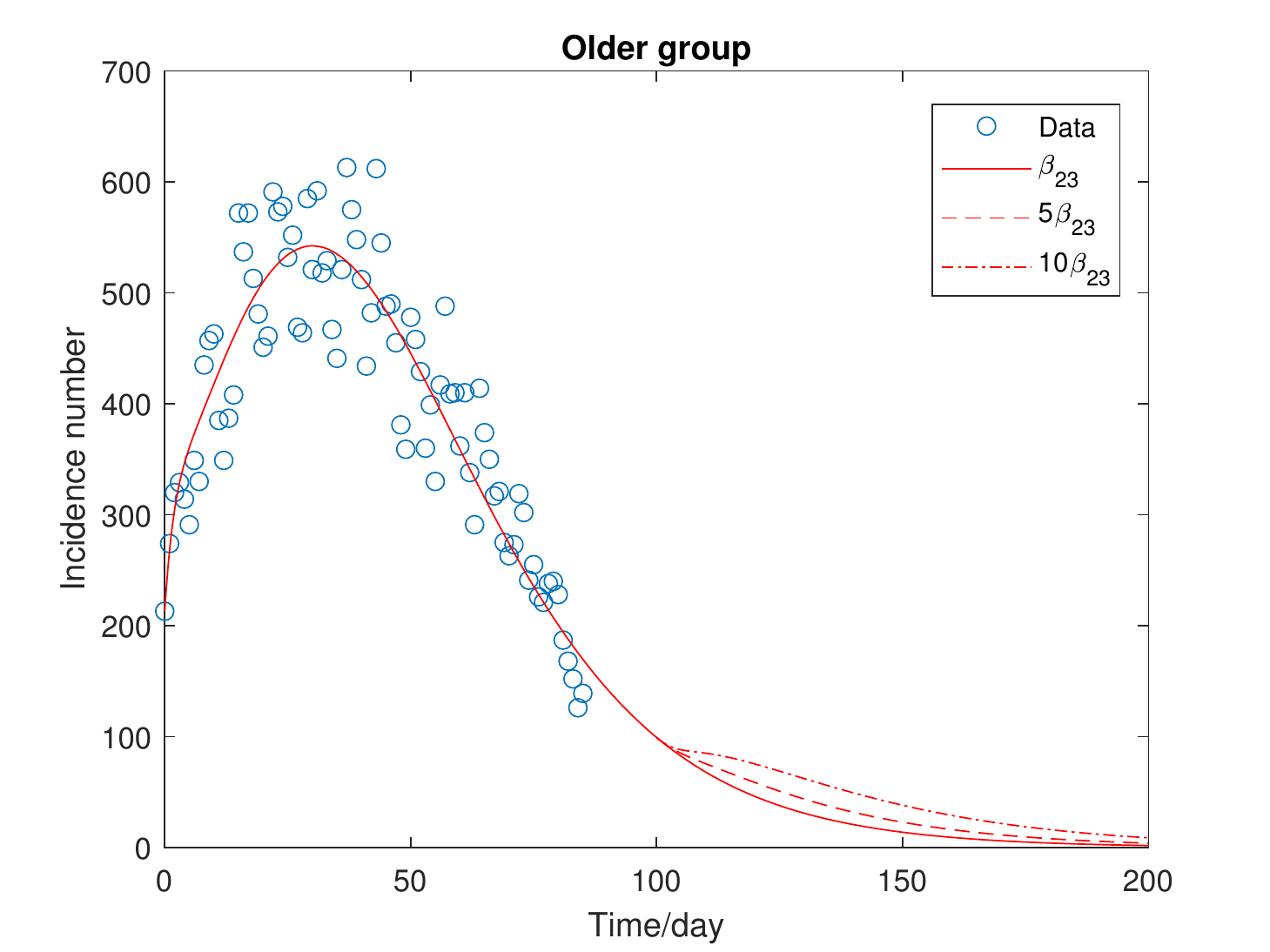}\label{fig:subfigure51}}\hfill
\subfloat[]{\includegraphics[width=0.5\textwidth]{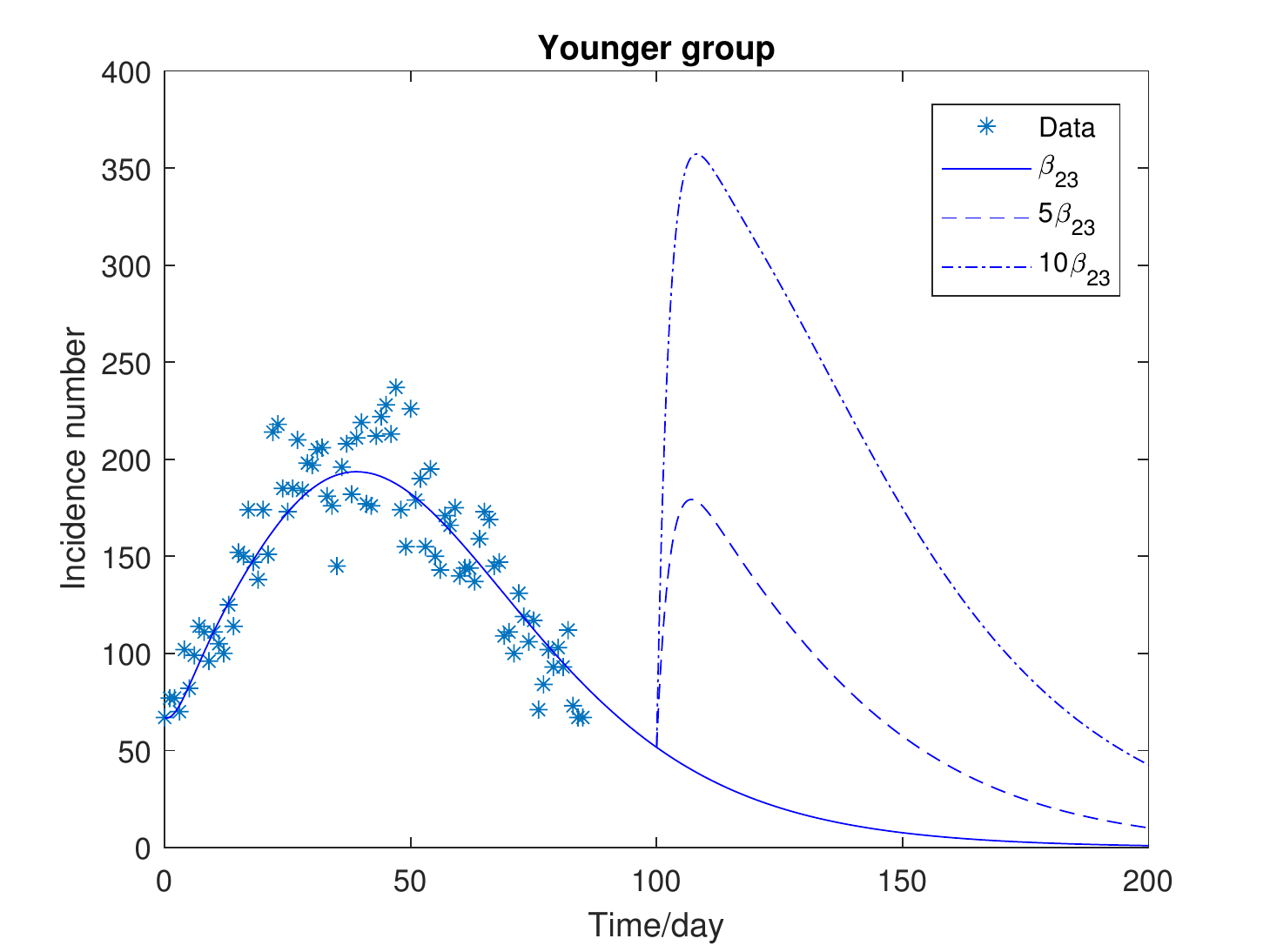}\label{fig:subfigure52}}\\
 \caption{Daily incidence number of the older group and the younger group respectively for different contact rate $\beta_{23}$. All the parameters remain unchanged as shown in \autoref{table:para} except $\beta_{23}.$}\label{fig:p5}
\end{figure}

We also examine how the between-group contact rate of the younger group $\beta_{23}$ may trigger the future epidemic waves. \autoref{fig:subfigure51} shows that the daily incidence number of the older group decreases monotonically if $\beta_{23}$ increases. However, by comparing \autoref{fig:subfigure41} and Figure \autoref{fig:subfigure51}, we observe that the daily incidence number of the older group is more sensitive to $\beta_{23}$ than $\beta_{21}$ even though a steady decline holds for both varying $\beta_{21}$ and $\beta_{23}$ in reasonable ranges. \autoref{fig:subfigure52} indicates that an increasing between-group contact rate of the younger group leads to a rapid increase in the daily incidence number of the younger group. The new wave reaches the peak of a larger size than the fourth wave in a short time and then declines.

Overall, \autoref{fig:p3}, \autoref{fig:p4}, \autoref{fig:p5} demonstrate that an increasing between-group contact rate in either the older group or the younger group may trigger future epidemic waves. The difference is that an increasing between-group contact rate of the older group leads to rapid growths of the daily incidence numbers of both groups whereas an increasing between-group contact rate of the younger group impacts more heavily on the younger group alone. The results indicate that a mixed transmission in different age groups plays an important role in triggering future waves and hence confirms the importance of vaccinating the younger group. The results also shed light on the demographic structure among the infected individuals in future waves: reported cases in the younger group will constitute a heavier proportion in the total infected individuals.

%P4

%P5

%P6

\noindent
(\rom{3}) Sensitivity analysis

In this section, we conduct sensitivity analysis to investigate the sensitivity of the total cumulative infection and infection number in each age-group on the estimated parameters. For each parameter, Latin Hypercube Sampling \cite{blower1994sensitivity,helton2006survey} is adopted to generate parameter values with assumed ranges and distributions as specified in \autoref{tab:sen}. We generate 3000 sets of parameter values for the analysis. By using these sets of values, partial rank correlation coefficients (PRCC) are calculated to determine the impact of varying parameters on the number of cumulative infections of COVID-19 between August and October, 2021 \cite{hamby1994review}. The PRCC indices range between -1 and 1, with positive (negative) values indicating a positive (negative) relationship and magnitudes indicating the relative level of impact on the quantity of interest, with a magnitude of 0 having almost no impact and 1 having the most influential impact.

\begin{table}[!htb] \caption {Parameter ranges and distributions for sensitivity analysis.} \label{tab:sen}
\centering
 \begin{threeparttable}
 \begin{tabular}{cccc}\toprule
  {Parameter} & {Lower Bound} & {Upper Bound} & {Distribution} \\ \midrule
%$\omega_2$ &0 & $^*\omega_1$ & U\\
  $\beta_{11}=\beta_{12}$ &0 & $10^{-5}$ & T\\
  $\beta_{13}=\beta_{14}$ &0 & $10^{-5}$ & T\\
  $\beta^v_{11}=\beta^v_{12}$ &0 & $10^{-5}$ & T\\
  $\beta^v_{13}=\beta^v_{14}$ &0 & $10^{-5}$ & T\\
  $\epsilon$ &0 & 1 & T\\
  $\beta_{21}=\beta_{22}$ &0 & $10^{-5}$ & T\\
  $\beta_{23}=\beta_{24}$ &0 & $10^{-5}$ & T\\
  $\omega$ &0 & $0.3$ & T\\
  $\tau_1$ &0 & 0.3 & T\\
  $\tau_2$ &0 & 0.3 & T\\
 \bottomrule
\end{tabular}
\begin{tablenotes}
      \small
      \item $T$ indicates triangular distribution with its peak value from \autoref{table:para}.
      \end{tablenotes}
\end{threeparttable}
\end{table}

\begin{figure}[!htb]
\captionsetup[subfigure]{labelformat=parens,labelfont=bf}
\centering
\subfloat[]{\includegraphics[width=0.5\textwidth]{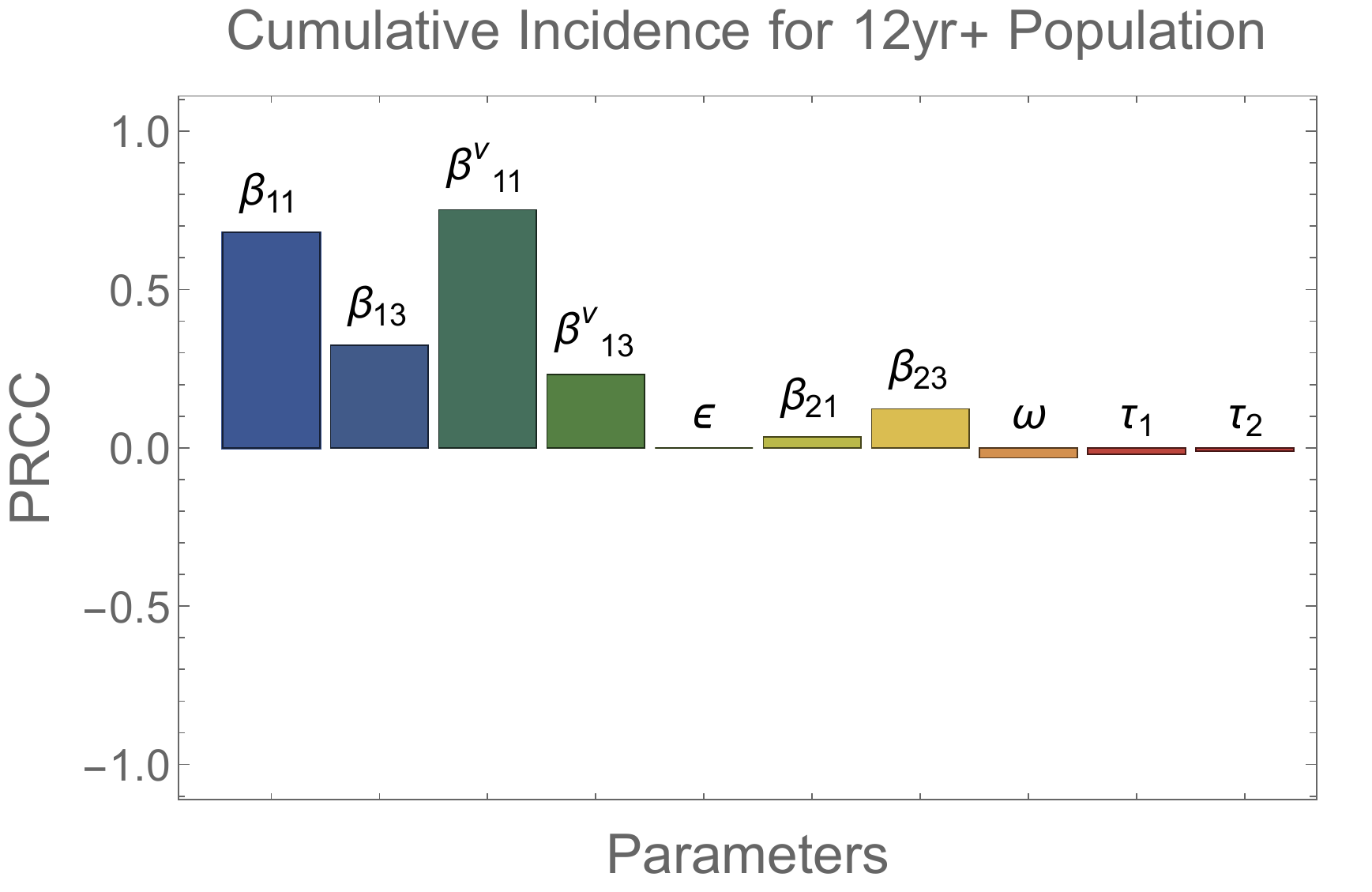}\label{fig:PRCC1}}\hfill
\subfloat[]{\includegraphics[width=0.5\textwidth]{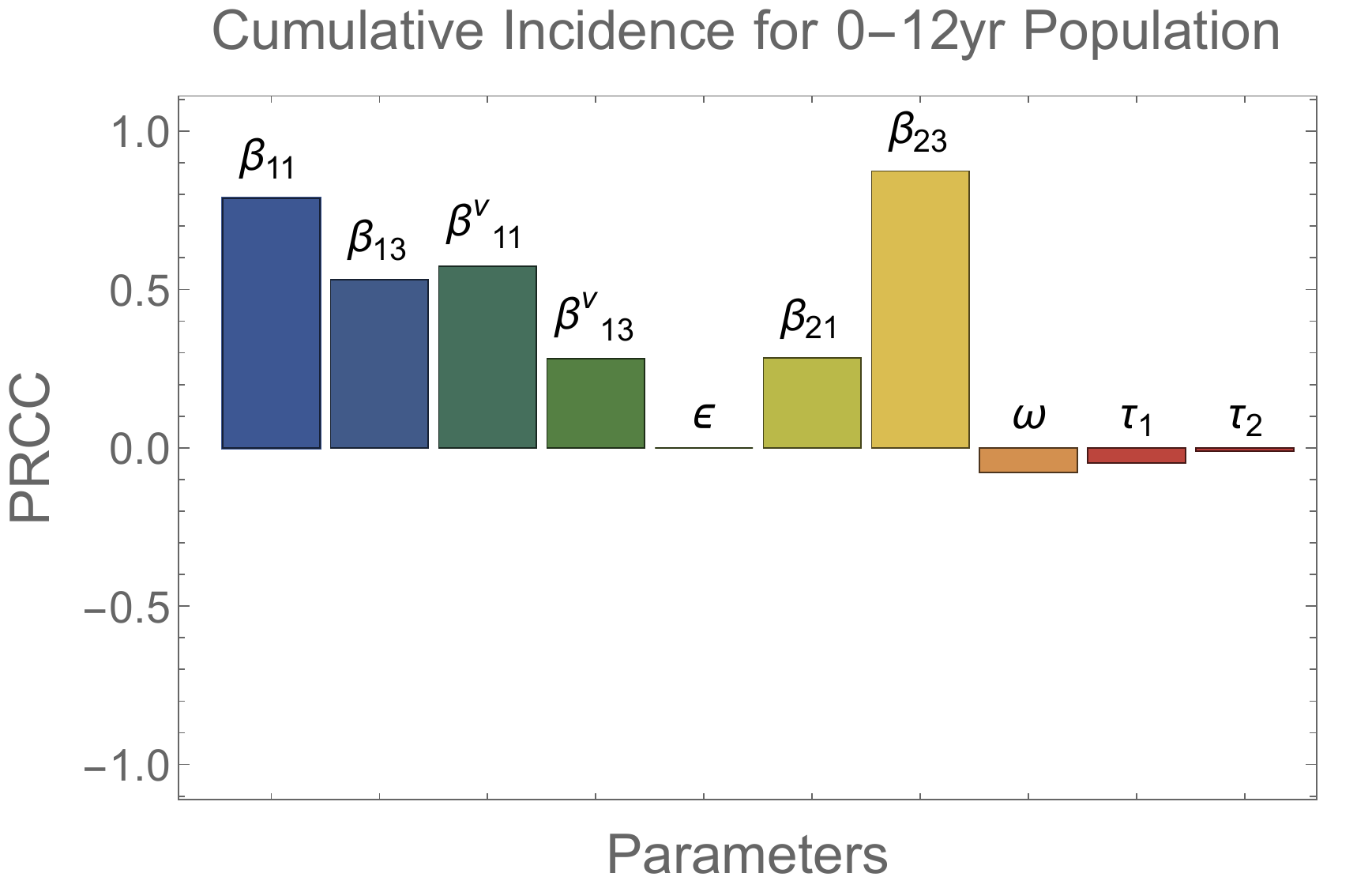}\label{fig:PRCC2}}\\
\subfloat[]{\includegraphics[width=0.5\textwidth]{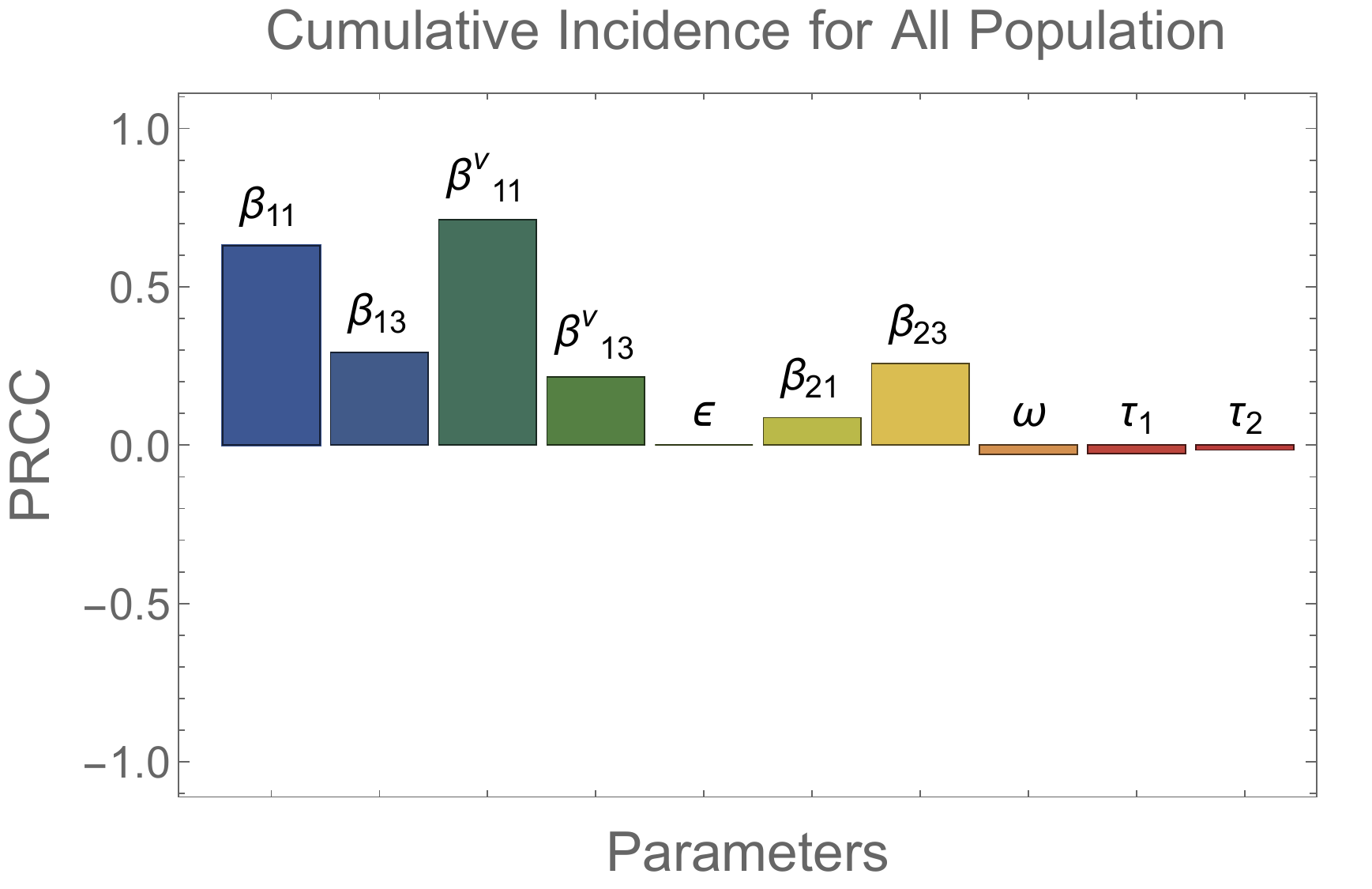}\label{fig:PRCC3}}\\
 \caption{Partial rank correlation coefficients (PRCC) calculated using parameter values from Latin Hypercube Sampling with respect to cumulative infections for population 12+ years of age, 0-12 years of age and all from Aug. 1, 2021 to Oct. 25, 2021 in Ontario, Canada.}\label{fig:PRCC}
\end{figure}

\autoref{fig:PRCC} shows the PRCC for cumulative incidence number of COVID-19 infection for the older group, the younger group, and the entire population in Ontario between August 1 and October 25, 2021. \autoref{fig:PRCC1}, \autoref{fig:PRCC2}, \autoref{fig:PRCC3} demonstrate that the transition rate $\tau_1, \tau_2$ from unreported infected class to reported infected class in respective group and vaccination rate $\omega$ in the older group have negative impact on the number of cumulative infections. The result is not surprising because individuals in the reported infected class are excluded from social activities due to self-isolation or hospitalization and therefore unable to transmit the pathogen whereas vaccinated individuals gain immunity and become less susceptible to infection. Comparing $\tau_1,\tau_2,\omega$, vaccination rate $\omega$ has a stronger influence than the reporting rates $\tau_1$ and $\tau_2$, which implies that vaccination is a more effective means to control the spread of disease than reporting infected cases. Moreover, reporting rate of the older group $\tau_1$ has a relatively larger impact than reporting rate of the younger group $\tau_2$ because the majority of infected population are adults.

\section{Conclusion and Discussion}

Since the initial identification of COVID-19 in Wuhan, China, the transmissible disease quickly escalates and becomes a global pandemic for about 2 years up to now. In Ontario, Canada, during the first year of the disease prevalence, the main intervention strategies have been nonpharmaceutical, such as provincial-wise lock down, keeping social distancing etc. due to the lack of effective vaccines. Such interventions are effective in mitigating the disease spread but cause significant economic loss at the same time and therefore are not long-lasting.

Progressing to the second year of the pandemic, in Ontario, starting from December 2020, the vaccine program gradually rolls out. The initial phase of the vaccine rollout targets the seniors population or vulnerable individuals and then gradually expands to individuals who are above the age of 12. The vaccine for children between the age of 5 to 11 was approved by Health Canada in late November but the rollout takes time and only $3.2\%$ of the children in the age of 5-11 are fully vaccinated by now \cite{ONdata}.

In this paper, we propose an age-stratified model that divides the entire population in Ontario into two groups: the older group where individuals are above the age of 12 and are eligible to receive vaccines and the younger group where individuals are below the age of 12 and are not eligible to receive vaccines. We fit the model to the daily incidence data of each group in Ontario between August 1, 2021 and October 25, 2021 and obtain a good fitting result.

The results demonstrate that between-group contact rate plays a more important role in triggering future waves than the within-group contact rate. The increasing between-group contact rate of the older group shifts the daily incidence number of either group from the decline to a rapid growth. However, an increasing between-group contact rate of the younger group mainly leads to a new epidemic wave in the younger group only. The results confirm the importance of achieving a high vaccine coverage in the younger group in order to mitigate the disease spread.

A new variant B.1.1.529 emerges in late November, 2021 in multiple African countries and quickly spreads to Canada in December. Early data indicate that the new variant is highly transmissible compared to other variants. Moreover, data show that currently, vaccines offer a much lower protection efficacy for even fully vaccinated individuals. Our model can be extended to study the invasion of the new variant, which leaves as future work.

\section*{Acknowledgement} XW is grateful for funding from the NSERC of Canada (RGPIN-2020-06825 and  DGECR-2020-00369).

%===================================================

%====================================================

%\bibliography{TwoGroupModel}{}
%\bibliographystyle{unsrtnat}

\end{document}